\documentclass[11pt]{article}
\usepackage{amsmath,amsfonts,theorem,epsfig}
\theoremstyle{plain}
\newtheorem{theorem}{Theorem}[section]
\newtheorem{lemma}[theorem]{Lemma}
\newtheorem{cor}[theorem]{Corollary}

\newcommand{\qed}[1]{\nopagebreak[4]\begin{flushright} 
\fbox{#1} \end{flushright}\pagebreak[2]}
\newcommand{\proof}{\noindent {\bf Proof}}
\newcommand{\hthree}{{\mathbb H}^3}

\newcommand{\bpq}{{\bf p},{\bf q}}

\newcommand{\cC}{{\cal C}}

\newcommand{\cH}{{\cal H}}
\newcommand{\cx}{{\mathbb C}}
\newcommand{\integers}{{\mathbb Z}}
\newcommand{\del}{\partial}
\newtheorem{Conjecture}[theorem]{Conjecture}
{\theorembodyfont{\rmfamily} 
\begin{document}
\title{Self-bumping of deformation spaces of hyperbolic 3-manifolds}
\author{K. Bromberg\footnote{Partially supported by a
grant from the Rackham School of Graduate Studies, University of
Michigan and by the Clay Mathematics Institute}\ \ and
J. Holt\footnote{Partially supported by a National Science Foundation
Postdoctoral Fellowship}}
\date{September 4, 2000}
\maketitle
\begin{abstract}
Let $N$ be a hyperbolic 3-manifold and $B$ a component of the interior
of $AH(\pi_1(N))$, the space of marked hyperbolic 3-manifolds homotopy
equivalent to $N$. We will give topological conditions on $N$
sufficient to give $\rho \in \overline{B}$ such that for every small
neighborhood $V$ of $\rho$, $V \cap B$ is disconnected. This implies
that $\overline{B}$ is not manifold with boundary.
\end{abstract}
\section{Introduction}
In this paper we study aspects of the topology of deformation spaces
of Kleinian groups.  The basic object of study is $AH(\pi_1(N))$, the
space of isometry classes of marked, complete hyperbolic 3-manifolds homotopy
equivalent to $N$, where $N$ is a compact, orientable, irreducible,
atoroidal 3-manifold with boundary.  The study of the global topology
of $AH(\pi_1(N))$ was begun by Anderson, Canary and McCullough in
\cite{ACM} for the case in which $N$ has
incompressible boundary.  They described necessary and sufficient
criteria for two components of the interior of $AH(\pi _1(N))$ to
``bump''; that is, to have intersecting closures.  We address the
question of when a component of the interior ``self-bumps''; that is,
if $B$ denotes such a component, then when is there an element $\rho$
in the closure of $B$ such that for any sufficiently small
neighborhood $V$ of $\rho$ in $AH(\pi _1(N))$ the set $V\cap B$ is
disconnected?  
In this paper we will establish the following result:

\medskip

\noindent
{\bf Theorem \ref{theoremA}} {\em
Let $N$ be a compact, orientable, atoroidal, irreducible 3-manifold
with boundary.  Suppose that $N$ contains an essential, boundary
incompressible annulus whose core curve is not
homotopic into a torus boundary component of $\del N$.  Let $B$ be a
component of the interior of $AH(\pi _1(N))$.  Then there is a
representation $\rho$ in $\overline{B}$ such that for any sufficiently
small neighborhood $V$ of $\rho$ in $AH(\pi _1(N))$ the set $V\cap B$
is disconnected.}

\medskip

Note that this result applies even when $N$ has compressible
boundary.  

In \cite{McMullen:complexeq} McMullen, using projective
structures and ideas of Anderson and Canary, proved Theorem \ref{theoremA}
when $N$ is an oriented $I$-bundle over a surface.  Our
techniques avoid the use of projective structures, and furthermore,
even in the $I$-bundle case we will find bumping representations that
are not detected with McMullen's methods. In a sequel, will we use the
techniques developed here to study the topology of the space of projective
structures with discrete holonomy.

We sketch the proof of Theorem \ref{theoremA} in the case where $N = S
\times [0,1]$ is an $I$-bundle over a closed surface of genus $\geq
2$. In this case the interior of $AH(\pi_1(N))$ consists of a single
component of quasifuchsian structures on $M= int\ N$, which is usually denoted
$QF(S)$.

\begin{figure}[ht]
\begin{center}
\epsfig{file=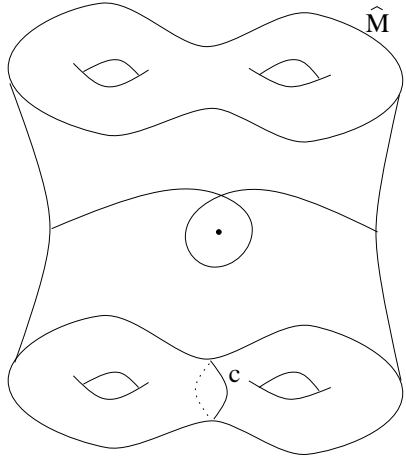,angle=-90,scale=.7}
\caption{The surface $S$ is immersed in $\hat{M}$ and is not homotopic
to and embedding.}
\label{surfacepic}
\end{center}
\end{figure}

To construct the representation where bumping occurs we start with a
hyperbolic structure on $M$ with a curve removed. That is choose a
simple closed curve $c$ on $S$ and let $\hat{M} = M - (c \times
\{1/2\})$. Then give $\hat{M}$ a geometrically finite hyperbolic
structure. Now, $\pi_1(\hat{M})$ has many conjugacy classes of
subgroups isomorphic to $\pi_1(S)$, for example $S \times \{1/4\}$ and
$S \times \{3/4\}$ each define such a subgroup. However, to find our
bumping representation we choose a non-standard subgroup of
$\pi_1(\hat{M})$ by wrapping $S$ around the removed curve (see Figure
\ref{surfacepic}).
Then the
hyperbolic structure on $\hat{M}$ defines a representation of
$\pi_1(\hat{M})$ and our choice of subgroup defines a representation
$\rho_\infty$ of $\pi_1(S)$. The cover $M_\infty$ associated to
this subgroup will be homeomorphic to $M$.

The next step is to construct an immersion $f:N \longrightarrow
\hat{M}_\infty$ in the
homotopy class associated to $\rho_\infty$ and then use $f$ to pull back a
hyperbolic structure $N_\infty$ on $N$. 
For each $\rho \in AH(\pi _1(N))$ there is a hyperbolic 3-manifold
$M_\rho$ homeomorphic to $M$. Given a small neighborhood $V$ of
$\rho_\infty$, for each $\rho \in V$ a general theorem allows us to
construct a smoothly varying family of hyperbolic structures
$N_\rho$ on the compact manifold $N$. Here $N_\rho$ has holonomy
$\rho$ and $N_{\rho_\infty} = N_\infty$. Since $N_\rho$ and $M_\rho$
have the same holonomy there 
will be an isometric immersion $f_\rho$ of $N_\rho$ in $M_\rho$. If
$\rho \in V \cap QF(S)$ then $c$ will have a geodesic representative
$c_\rho$ in $M_\rho$ and there will be a canonical homeomorphism
between $M_\rho - c_\rho$ and $\hat{M}$. Furthermore, geometric
considerations will show that the image of $f_\rho$ misses $c_\rho$ so
we can view $f_\rho$ as a map to $\hat{M}$. In particular, we can
compare the homotopy classes of the maps $f_\rho$ in $\hat{M}$.

The heart of the proof is that we can find representations $\rho_0$
and $\rho_1$ in $V \cap QF(S)$, such that $f_0$ is homotopic in
$\hat{M}$ to the original immersion $f$ while $f_1$ is homotopic to
an embedding. If $\rho_0$ and $\rho_1$ are in the
same component of $V \cap QF(S)$ then our smoothly varying family of
hyperbolic structures $N_\rho$ will define a homotopy between $f_0$
and $f_1$ in $\hat{M}$. This contradiction proves the theorem.

To find the representation $\rho_0$ we take a small deformation of
$\hat{M}$ that fills in the torus boundary to give a manifold
homeomorphic to $M$. To find $\rho_1$ we take a small deformation of
$M_\infty$ that resolves the rank one cusp. In this case the
homeomorphism type is preserved. Since $M_0$ is geometrically very
close to $\hat{M}$, $f_0$ will have the same homotopy class as $f$. In
$M_\infty$, $f$ lifts to an embedding and therefore $f_1$ will be an
embedding in $M_1$.

It is worthwhile to compare this result with the bumping of distinct
components examined in \cite{Anderson/Canary} and
\cite{ACM}. As mentioned above in
\cite{ACM},
necessary and sufficient conditions are given for components to
bump. We will not state them here, but at the very least we need a
manifold with more topology than an $I$-bundle so that the interior of
$AH(\pi_1(N))$ will have more than one component. The construction of
the bumping representation is then very similar to the one above. 

We first remove a suitably chosen simple closed curve $c$ from $M =
int\ N$ to obtain an new manifold $\hat{M}$.  We then find a cover
$M_\infty$ of $\hat{M}$ that is homotopy equivalent, but in this case
not homeomorphic to, $M$. A hyperbolic structure on $\hat{M}$ defines a
hyperbolic structure on 
$M_\infty$. As above we make a small deformation $M_0$ of $\hat{M}$
that will be homeomorphic to $M$ while a small deformation $M_1$ of
$M_\infty$ will be homeomorphic to $M_\infty$.  Although $M_0$ and
$M_1$ are not homeomorphic, their holonomy representations $\rho_0$ 
and $\rho_1$ will both be near the holonomy representation $\rho_\infty$ of
$M_\infty$. The next, and last, step is the real difference between
the two arguments. As the components of the interior of $AH(\pi_1(N))$ are
parameterized by the (marked) oriented homeomorphism types of $N$, $\rho_0$ and
$\rho_1$ must be in distinct components that bump at $\rho_\infty$. 

\noindent{\bf Acknowledgments.}

The authors would like to thank Jeff Brock and Dick Canary for
interesting and helpful discussions.

\vspace{3mm}

\section{Preliminaries}\label{Preliminaries}

A {\it Kleinian group} is a discrete, torsion free subgroup of the 
orientation preserving isometries of hyperbolic 3-space, $\hthree$.
In the upper-half-space model of $\hthree$ the orientation-preserving
isometries are identified with the group $PSL_2(\mathbb C)$, so that a
Kleinian group can be considered a discrete, torsion free subgroup of
$PSL_2(\mathbb C)$.   

Let $\Gamma$ be a Kleinian group and set $M$ to be the quotient
manifold $\hthree/\Gamma$.  The {\it convex core} of $M$ is the
smallest convex submanifold of $M$ whose inclusion in $M$ is a
homotopy equivalence.  If the convex core has finite volume, and
$\Gamma$ is finitely generated then $\Gamma$ is called {\it
geometrically finite}.  In addition, a geometrically finite Kleinian
group is {\it minimally parabolic} if every maximal parabolic subgroup
is of rank 2. 

Let $R(\pi_1(N)) = Hom(\pi_1(N), PSL_2(\cx))/PSL_2(\cx)$ be the
space of conjugacy classes of representations of $\pi_1(N)$ in
$PSL_2(\cx)$ where $N$ is a compact, orientable, atoroidal
3-manifold. The subset $AH(\pi_1(N)) \subset R(\pi_1(N))$ consists of
the discrete, 
faithful representations of $\pi_1(N)$, modulo conjugacy.  It is a
result of J\o rgensen \cite{Jorgensen:AHclosed} that $AH(\pi _1(N))$ is
a closed subset of $R(\pi _1(N))$.  By work of Marden
\cite{Marden} and Sullivan \cite{Sullivan} the interior of $AH(\pi
_1(N))$ is $MP(\pi _1(N))$, the minimally parabolic representations. 

A representation $\rho \in AH(\pi_1(N))$ determines an oriented
hyperbolic 
manifold $M_\rho = \hthree/\rho(\pi_1(N))$ along with a homotopy
equivalence, $f_\rho :N \longrightarrow M_\rho$. While in general
$MP(\pi_1(N))$ will have many components, in
this paper our interest is the topology of the closure of a single
component $B$. Note that $AH(\pi_1(N))$ is determined only by the
homotopy type of $N$. We can therefore assume that $N$ is chosen such
that if $\rho$ is in $B$ then there is a homeomorphism from $M_\rho$
to the interior of $N$ that is a homotopy inverse for
$f_\rho$. We can also orient $N$ such that this homeomorphism is
orientation preserving. Then $B$ will be the unique component of
$MP(\pi_1(N))$ satisfying these two properties.

We also need to work with hyperbolic structures on the compact
manifold $N$ that may not extend to complete hyperbolic structures on
an open manifold containing $N$. We let $\cH(N)$ be the space of
hyperbolic metrics
on $N$. Given two hyperbolic metrics on $N$ the identity map will be a
biLipschitz map between the two metrics. Given a structure, $N' \in
\cH(N)$, a neighborhood $N'(\epsilon)$ of $N'$ consists of those
structures in $\cH(N)$ for which the identity map from $N'$ is a $(1 +
\epsilon)$-biLipschitz map. The $N'(\epsilon)$ are a basis of
neighborhoods for $N'$.

Theorem 1.7.1 in \cite{Canary:Epstein:Green} describes the local
structure of a neighborhood of $N'$. We will need the following simple
consequence of this theorem:

\begin{theorem}{\cite{Canary:Epstein:Green}}
\label{localproduct}
The holonomy map $\cH(N) \longrightarrow R(\pi_1(N))$ is locally
onto. Furthermore, for any neighborhood $V$ of $N'$ there exists a
neighborhood $U \subset V$, such that if $N_0$ and $N_1$ are hyperbolic
structures in $U$ with holonomy $\rho_0$ and $\rho_1$, respectively,
and $\rho_t$, $0 \leq t \leq 1$, is a path in the image of $U$ then
there is a path $N_t$ in $U$, where each $N_t$ has holonomy $\rho_t$.
\end{theorem}

Now assume that $\del N$ contains at least one torus, $T$.  Choose a
meridian and longitude for this 
torus  such that elements of $\pi_1(T) = \integers \oplus \integers$
are determined by a pair of integers.  Let $(p,q)$ be a pair of
relatively prime integers.  Let 
$N(p,q)$ denote the result of performing $(p,q)$-Dehn filling on $N$
along this torus; that is, there exists an embedding $d_{p,q}: N
\longrightarrow N(p,q)$ such that $\overline {N(p,q)-d_{p,q}(N)}$ is a solid torus bounded
by $d_{p,q}(T)$ and the image of the $(p,q)$ curve on $T$ is trivial
in $N(p,q)$. Let $\gamma$ denote the core curve of of the solid
torus.  If $N$ and $N(p,q)$ have complete hyperbolic structures, 
$M$ and $M(p,q)$, on their interiors then $M(p,q)$ is a {\em hyperbolic
Dehn filling} of $M$ if $M(p,q) - d_{p,q}(M)$ contains the geodesic
representative of $\gamma$. Note that a hyperbolic structure $M(p,q)$
may not be a hyperbolic Dehn filling of $M$ if $\gamma$ is not
isotopic to its geodesic representative. Also note that
the holonomy representation $\rho$ for $M(p,q)$ induces a
non-faithful, holonomy representation, $\rho_{p,q}$, for $N$ via
pre-composition with $(d_{p,q})_*$. 

If $N$ has $k$ torus boundary components, we can Dehn fill each of
them. Let relatively prime integers, $(p_i, q_i)$, be the Dehn filling
coefficients for the $i$-th torus and let $({\bf p}, {\bf q}) =
(p_1,q_1; \dots; p_k,q_k)$. Then $N({\bf p}, {\bf q})$ is the $({\bf
p},{\bf q})$-Dehn filling of $N$.

The following theorem has an extensive history.  The interested reader
should also see \cite{Benedetti:Petronio}, \cite{Thurston:Notes},
\cite{Bonahon:Otal:arbshort}, and \cite{Comar:thesis}.

\begin{theorem}{The Hyperbolic Dehn Surgery Theorem
(\cite{Bromberg:dehn})}
\label{surgery} 

Let $M$ be a compact 3-manifold with $k$ torus boundary components and
assume $M$ has a minimally parabolic hyperbolic structure with
holonomy $\rho$. We then have the following:
\begin{enumerate}
\item Except for a finite number of 
pairs for each $i=1,\ldots, k$, for each collection of relatively prime pairs
$(\bpq)$ there exist a geometrically finite hyperbolic $(\bpq)$-Dehn
filling $M(\bpq)$ of $M$.

\item $\rho_{\bpq} \rightarrow \rho$ as $|\bpq| \rightarrow \infty$
$(|\bpq| = |p_1| + |q_1| + \cdots + |p_k| + |q_k|).$ 

\item If $X$ is the complement of a neighborhood of the cusps and $|\bpq| >
n$ then $d_{\bpq}|_{X}$ is $K_n$-biLipschitz with $K_n \rightarrow 1$
as $n \rightarrow \infty$.

\end{enumerate}
\end{theorem}

\section{Wraps and twists}
\label{wraptwist}

Let
$$X = [-1,1] \times [-1,1] \times S^1$$
and
$$\hat{X} = X - ([-\frac {1}{3},\frac {1}{3}] \times [-\frac {1}{3},\frac {1}{3}] \times S^1).$$
We begin be defining maps of the annulus,
$$A= [-1,1] \times S^1$$
into $\hat{X} \subset X$. First we define $w:A \longrightarrow \hat{X}$
by 
$$w(x,\theta) = (-\frac{1}{2}\sin(\pi x), \frac{1}{2}\cos(\pi x),
\theta).$$
\begin{figure}[ht]
\begin{center}
\epsfig{file=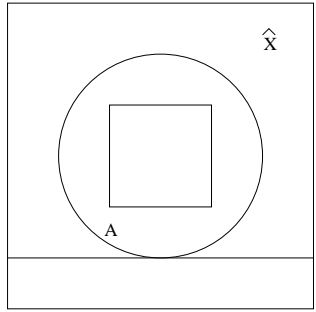,angle=-90,scale=.7}
\caption{The image of $A$ under the map $w_1$ in a cross section of
$\hat{X}$.}
\label{wrappic}
\end{center}
\end{figure}
We next define a sequence of maps $w_n:A \longrightarrow
\hat{X}$ for each $n > 0$. For each $t$ and $t'$ with $-1\leq
t<t'\leq 1$ we let $h_{t,t'}:([t,t'] \times S^1) \longrightarrow A$ be
a homeomorphism that satisfies the conditions, $h_{t,t'}(t,\theta) =
(-1,\theta)$ and $h_{t,t'}(t',\theta) = (1,\theta)$. To define $w_n$
we choose real numbers, $t_0, \dots, t_n$ with $-\frac {1}{3}=t_0 < t_1 <
\cdots <t_n=\frac {1}{3}$, and let
\[w_n(x,\theta) = \left\{ \begin{array}{ll}
\left(\frac{3}{2}x + \frac12,-\frac12, \theta\right) & \mbox{if $-1 \leq x <
-\frac13$} \\
w \circ h_{t_i,t_{i+1}} & \mbox{if $t_i \leq x < t_{i+1}$}\\
\left(\frac32 x- \frac12,-\frac12, \theta\right) & \mbox{if $\frac13
\leq x \leq 1$.}
\end{array}
\right. \]
The map $w_n$ wraps the annulus $n$ times around the missing core of
$\hat{X}$. For $n=0$, we define $w_0$ by $w_0(x, \theta) =
(x,-1/2,\theta)$.

Our next family of maps, $t_{n,m}:\hat{X} \longrightarrow \hat{X}$,
are homeomorphisms which {\em Dehn twist} $\hat{X}$. They are defined
by the following formula:
\[ t_{n,m} = \left\{ \begin{array}{ll}
(x,y,\theta) & \mbox{if $-1 \leq x < -\frac {1}{3}$ or $\frac {1}{3} < x \leq 1$}\\
(x,y,\theta + 3n\pi(x+\frac {1}{3}) & \mbox{if $-\frac {1}{3} \leq x \leq \frac {1}{3}$ and
$y>\frac {1}{3}$}\\
(x,y,\theta + 3m\pi(x+\frac {1}{3}) & \mbox{if $-\frac {1}{3} \leq x \leq \frac {1}{3}$ and
$y<-\frac {1}{3}$.}
\end{array}
\right. \]

\begin{lemma}
\label{thetwist}
The maps $w_n$ and $t_{k(n+1),kn} \circ w_n$ are homotopic rel $\del
A$ for any positive integer $n$ and any integer $k$.
\end{lemma}

\proof

\begin{figure}[ht]
\begin{center}
\epsfig{file=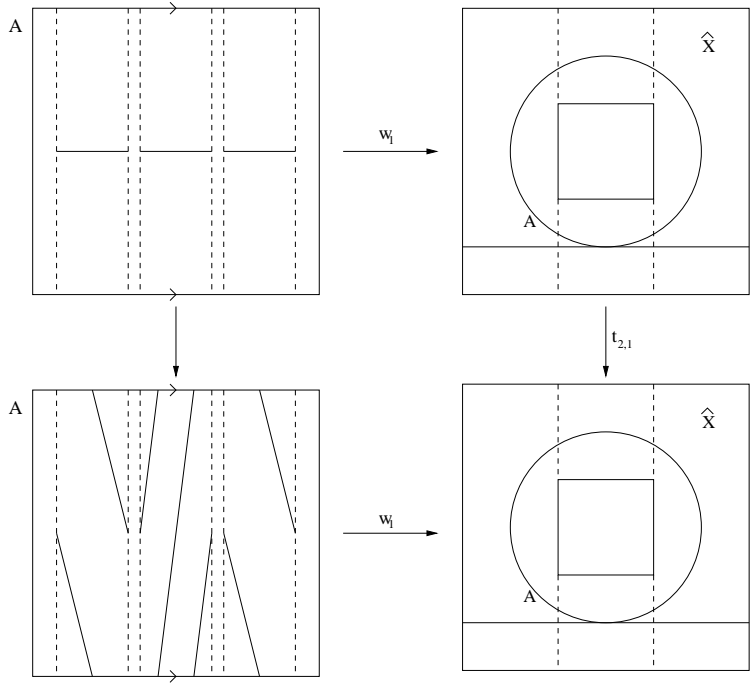,angle=-90,scale=.5}
\caption{By identifying the top and bottom of the squares on the left
we obtain (two copies of) the annulus $A$. The preimage of the
$w_1(A)\cap \hat X_{\frac {1}{3}}$ is the three dashed annuli. The effect
of $t_{2,1}$ on $A$, is two dehn twists on the center annuli and a
single dehn twist in the opposite direction on the two outside
annuli. As we see from the picture in the lower left, the net effect
on $A$ is a map that is homotopic to the identity.}
\label{dehntwist}
\end{center}
\end{figure}
Let $\hat X_{\frac {1}{3}}=([-\frac{1}{3}, \frac {1}{3}]\times [-1,
1]\times S^1)\cap \hat X$ denote the middle-third of $\hat X$; it has
two components, the upper half and the lower half.  The image of $A$ under the map $w_n$ intersects $\hat X_{\frac {1}{3}}$, so
that $w_n ^{-1}(w_n(A)\cap \hat X_{\frac{1}{3}})$ consists of $2n+1$
essential sub-annuli of $A$; $n$ of the
annuli map into the upper half of the middle third, while $n+1$ of the
annuli map into the lower half. On the each of the $n+1$ annuli
mapping into the lower half, $t_{k(n+1),kn}$ is a $kn$-Dehn twist,
while on the $n$ upper annuli $t_{k(n+1),kn}$ is a $-k(n+1)$-Dehn
twist. Therefore the total affect of $t_{k(n+1),kn}$ is a $kn(n+1) -
k(n+1)n=0$-Dehn twist and $w_n$ is homotopic to $w_n \circ
t_{k(n+1),kn}$ rel $\del A$ (see Figure \ref{dehntwist}).
\qed{proof of Lemma \ref{thetwist}}

We now relate the maps $t_{n,m}$ to the Dehn filling of
$\hat{X}$.  As
our coordinates for Dehn filling we choose the meridian to be the
unique homotopy class that is trivial in $X$ and the longitude to be the
curve $\{\frac {1}{3}\} \times \{\frac {1}{3}\} \times S^1$. Recall
the Dehn filling maps $d_{1,k}:\hat{X} \longrightarrow \hat{X}(1,k)$.

\begin{lemma}
\label{dehnit}
For each $t_{n,m}$ there exists a homeomorphism
$h_{n,m}:\hat{X}(1,n-m) \longrightarrow \hat{X}(1,0)$ such that
$d_{1,0} \circ t_{n,m} = h_{n,m} \circ d_{1,n-m}$.
\end{lemma}
\nopagebreak
\proof

The map $t_{n,m}$ takes the $(1,n-m)$-curve to the
$(1,0)$-curve so $d_{1,0} \circ t_{n,m}$ takes the $(1,n-m)$ to a
trivial curve in $\hat{X}(1,0)$. On the image of $\hat{X}$ in
$\hat{X}$, we define 
$h_{n,m}$ to satisfy the equation, $d_{1,0} \circ t_{n,m} = h_{n,m}
\circ d_{1,n-m}$. Since the image of the $(1,n-m)$ curve is trivial in
$\hat{X}(1,n-m)$, $h_{n,m}$ extends to a homeomorphism.
\qed{proof of Lemma \ref{dehnit}}

Let 
$$\del _0 X = [-1,1] \times \{-1,1\} \times S^1 \subset
X$$
and
$$\del_1 X = \{1\} \times [-1,1] \times S^1 \subset X.$$
Also assume that $N$ is a compact manifold with boundary and that $N$
contains an essential, boundary incompressible annulus. Then there is
a pairwise embedding of $(X, \del _0 X)$ in $(N, \del N)$ such that
$\del_1 X$ is an essential, boundary incompressible annulus.
Identify $A$ with the lower half of $\del_0 X$; that is, the 
annulus $[-1,1] \times \{-1\} \times S^1$. Let $c = \{0\} \times \{0\}
\times S^1$ be the core curve of $X$ and let $\hat{M} = M - c$ where
$M$ is the interior of $N$.

For each integer $n \geq 0$ we define an immersion $s_n:N
\longrightarrow M \subset N$ as follows. The map $s_n$ is homotopic to the
identity map and a homeomorphism onto its image outside of $X$. We
also require that $s_n(N) \cap c = \emptyset$ and that $s_n$ restricted
to $A$ is homotopic to $w_n$ rel $\del A$. This completely defines
$s_n$ up to homotopy in $\hat M$.  We call any map that satisfies
these properties a {\em shuffle immersion}.

\begin{figure}[ht]
\begin{center}
\epsfig{file=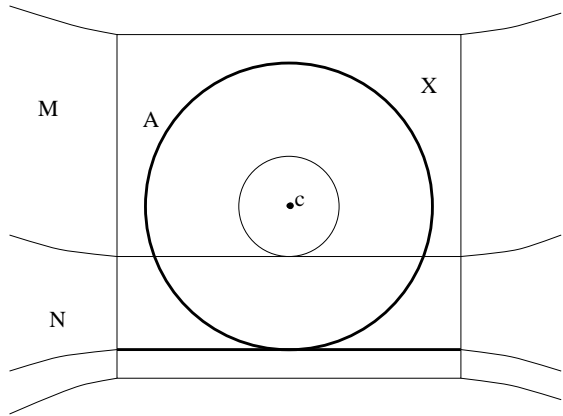, angle=-90,scale=.6}
\caption{The map $s_1$ immerses $N$ in $M$ and is not homotopic to an
embedding in $\hat{M}$.}
\end{center}
\end{figure}

\begin{lemma}
\label{immersion}
A shuffle immersion $s_n$ satisfies the following properties:
\begin{enumerate}
\item If $n \neq m$ then $s_n$ and $s_m$ are not homotopic in $\hat{M}$.

\item For each integer $k$, there is an orienation preserving
homeomorphism $h_k 
:\hat{M}(1,k) \longrightarrow M$ such that $s_n$ and $h_k \circ
d_{1,k}\circ s_n$ are homotopic in $\hat{M}$. Here, $M = \hat{M}(1,0)$.

\item The cover of $\hat{M}$ associated to $(s_n)_*(\pi_1(N))$ is
homeomorphic to $M$ and $s_n$ lifts to an embedding $\hat
s_n:N\longrightarrow M$ which is homotopic to $s_0$ in $\hat{M}$.
\end{enumerate}
\end{lemma}

\proof
\begin{enumerate}
\item If $n \neq m$, the maps $(s_n)_*(\pi_1(N))$ and
$(s_m)_*(\pi_1(N))$ are non-conjugate subgroups of
$\pi_1(\hat{M})$  and therefore the maps $s_n$ and $s_m$ are not
homotopic.

\item On $\hat{X}(1,k) \subset \hat{M}(1,k)$ we let $h_k =
h_{k(n+1),kn}$. Using Lemma \ref{dehnit}, we see that $h_k$ extends to
a homeomorphism from $\hat{M}(1,k)$ to $M$. By Lemma \ref{thetwist},
$s_n$ and $h_k \circ d_{1,k}\circ s_n$ are homotopic
in $\hat{M}$.

\item This is an easy exercise in 3-manifold topology which we leave
to the reader.
\end{enumerate}
\qed{proof of Lemma \ref{immersion}}

\section{Self-bumping}
\label{selftouching}
We now use the topology we developed in \S \ref{wraptwist}. With the
same assumptions as in \S \ref{wraptwist} we fix a shuffle immersion
$f=s_d$, with $d >0$. Note that such a shuffle immersion exists if and
only if $N$ contains an essential, boundary incompressible
annulus. However, for $M$ and $\hat{M}$ to support complete
hyperbolic structures we need to make further topological
restrictions. Namely, $N$ must be irreducible and atoroidal and the
simple closed curve $c$ must be primitive and not homotopic to a
torus boundary component of $\del N$. Then $M$ and $\hat{M}$ satisfies the
conditions of Thurston's hyperbolization theorem (see Lemma 2.5.10
in \cite{Holt:thesis}) and we fix a minimally parabolic hyperbolic
structure $\hat{M}_\infty$ on $\hat{M}$ with holonomy representation
$\hat{\rho}_\infty$. We also let $N_\infty$ be the hyperbolic metric on
$N$ obtained as the pull-back by $f$ of the metric $\hat{M}_\infty$ on
$\hat{M}$.  

We now set up a notational system that will hold for the remainder of
the paper.  For an index $\alpha$, $N_\alpha$ is a
hyperbolic structure on $N$ and $\rho_\alpha$ will be the associated
holonomy representation. Let $M$ be the interior of $N$. As we noted
in the introduction, if $\rho_\alpha \in AH(\pi_1(N))$ then $M_\alpha$ is a complete
hyperbolic structure, marked by $N$.
As $N_\alpha$ has the same holonomy as $M_\alpha$ there will be an
an isometric immersion, $f_\alpha:N_\alpha \longrightarrow M_\alpha$,
with $f_\alpha$ a 
homotopy equivalence. In other words, $f_\alpha$ is the marking
map. Let $c_\alpha$ denote the geodesic representative of $c$ in
$M_\alpha$. 

\begin{lemma}
\label{short}
Let $V$ be a small neighborhood of $\rho_\infty$. Then for each
$N_\alpha$ near $N_\infty$ with $\rho_\alpha \in V \cap MP(\pi_1(N))$,
$f_\alpha(N_\alpha) \cap c_\alpha = \emptyset$. 
\end{lemma}

\proof

By compactness there exists a K such that for any $p \in N$ we can find a
non-trivial simple closed curve $\gamma_p$ through $p$, and not
homotopic to $c$, 
with length $< K$ in $N_\infty$. We choose $V$ small enough such that
all structures in the neighborhood are 2-biLipschitz from
$N_\infty$. The Margulis lemma implies that there exists an $\epsilon$
such that, for any complete hyperbolic 3-manifold, if a homotopically
non-trivial simple closed curve intersects a homotopically distinct
geodesic of length $< \epsilon$ it has length $> 3K$.
Furthermore, since the length of curves is continuous on $R(\pi_1(N))$,
we can further shrink $V$ so that the curve, $c_\alpha$, has length $<
\epsilon$ and therefore $f_\alpha(\gamma_p)$, which has length $< 2K$,
does not intersect $c_\alpha$; implying that $p \not\in c_\alpha$.

\qed{proof of Lemma \ref{short}}

Recall that $B$ is a component of $MP(\pi _1(N))$ so that
the  marking map $f_\alpha$ has a homotopy inverse which is an
orientation preserving homeomorphism between $M_\alpha$ and $M$ if and
only if $\rho_\alpha \in B$.

\begin{lemma}\label{dehn}
For the shuffle immersion $f$, there exists a sequence of hyperbolic
structures $N_k$ with holonomy representations $\rho_k$, such that:
\begin{enumerate}
\item $N_k \rightarrow N_\infty$.

\item $\rho_k \rightarrow \rho_\infty$.

\item There exist homeomorphisms $h_k:M_k \longrightarrow M$ such
that $h_k$ is a homotopy inverse for $f_k|_M$, $h_k(c_k) = c$ and $f$
and $h_k \circ f_k$ are homotopic in $\hat{M}$. In particular, $\rho_k
\in B$.
\end{enumerate}
\end{lemma}

\proof

\begin{enumerate}
\item For large $n$,
let $M_n = \hat{M}_\infty (1,n)$ be the manifolds obtained by
performing hyperbolic Dehn surgery on $\hat{M}_\infty$ as in Theorem
\ref{surgery}. Since $f_\infty(N_\infty)$ is contained in a compact
subset of $\hat{M}$, Theorem \ref{surgery} also implies that the maps
$d_{1,n}:\hat{M}_\infty \longrightarrow M_n$ restricted to
$f_\infty(N_\infty)$ are $K_n$-quasi-isometries with $K_n \rightarrow
1$ as $n \rightarrow \infty$. The hyperbolic structures $N_n$
defined by pulling back the hyperbolic metric on $M_n$ by $d_{1,n} \circ
f_\infty$ converge to $N_\infty$.\\

\item Again, as in Theorem \ref{surgery}, the maps $d_{1,n}$ define holonomy
representations $\hat{\rho}_{1,n}$ of $\pi_1(\hat{M})$ with
$\hat{\rho}_{1,n} \rightarrow \hat{\rho}_\infty$. The holonomy
representations, $\rho_n$, are restrictions of $\hat{\rho}_{1,n}$ so
$\rho_n \rightarrow \rho_\infty$.\\

\item These homeomorphisms are supplied by Lemma \ref{immersion}.\\
\end{enumerate}
\qed{proof of Lemma \ref{dehn}}

The following lemma will be  used to detect when two representations
are not contained in the same component of $V \cap B$.

\begin{lemma}\label{homotopy}
Let $U$ be a neighborhood of $N_\infty$ that satisfies the conclusion
of Theorem \ref{localproduct} and Lemma \ref{short} and let $V$ be the
image of $U$ under the holonomy map. Let $N_0$ and $N_1$ be hyperbolic
structures in $U$ with holonomy $\rho_0$ and $\rho_1$, both in $V \cap
MP(\pi _1(N))$. Also assume that $h_i:M_i \longrightarrow M$, $i=0,1$, are
homeomorphisms that are homotopy inverses of $f_i|_M$ and $h_i(c_i) =
c$. If $\rho_t$, $0\leq t \leq 1$, is a path in $V \cap MP(\pi_1(M))$
then $h_0 \circ f_0$ and $h_1 \circ f_1$ are homotopic in $\hat{M}$.
\end{lemma}

\proof

By Theorem \ref{localproduct} we have a path of structures $N_t$ in
$V$ with holonomy $\rho_t$. The $\rho_t$ are in the same
component of $MP(\pi_1(M))$ as $\rho_0$ and $\rho_1$ are in, so there are
homeomorphisms, $h_t:M_t \longrightarrow M$, that are homotopy
inverses of $f_t|_M$. We can assume that the push-forward of the
hyperbolic metrics on $M_t$ to $M$ is a continuously changing family
of metrics on $M$. Furthermore, as all the $c_t$ are short geodesics,
they will be simple. Hence $h_t(c_t)$ is an isotopy of $c$ in
$M$. We can therefore modify the $h_t$ such that $h_t(c_t) = c$.
Then $h_t \circ f_t$ will vary continuously in $t$. By Lemma
\ref{short} $f_t(N_t) \cap c_t = \emptyset$ so $h_t \circ f_t$ is a
homotopy between $h_0 \circ f_0$ and $h_1 \circ f_1$ in $\hat{M}$.
\qed{proof of Lemma \ref{homotopy}}

We next apply Lemma \ref{homotopy} to show that distinct shuffle
immersion force $V \cap B$ to be disconnected.

\begin{lemma}\label{neighbourhood}
Let $f,f':N \longrightarrow \hat{M}\subset M$,  be distinct shuffle
immersions. Assume that there exists
minimally parabolic structures $\hat{M}_\infty$ and
$\hat{M}'_\infty$ on $\hat{M}$ such that the pulled-back hyperbolic
structures $N_\infty$ and $N'_\infty$ are isometric and hence define
the same holonomy representation, $\rho_\infty$. Then for every small
neighborhood $V$ of $\rho_\infty$, $V \cap B$ is disconnected.
\end{lemma}

\proof

Let $M_n$, $N_n$, $f_n$, $h_n$, and $\rho_n$ and $M_n'$, $N'_n$,
$f'_n$, $h'_n$, and
$\rho'_n$ be the hyperbolic structures, isometric immersions and
holonomy representations given by Lemma \ref{dehn} for $f$ and $f'$,
respectively. Choose an open neighborhood $V$ of $\rho_\infty$ given by
Lemma \ref{short}. 

There exists integers $n$ and $m$ such that $\rho_n, \rho'_m \in V$.
The intersection $V \cap B$ is an open subset of the manifold $B$ so
the connected components of $V \cap B$ are path 
connected.  If $\rho_n$ and $\rho'_m$ are in the same component of $V
\cap B$ then Lemma \ref{homotopy} implies that $h_n \circ f_n$ and
$h'_m \circ f'_m$ are homotopic in $\hat{M}$. On the other hand, by
Lemma \ref{dehn}, $h_n \circ f_n$ and $h'_m \circ f'_m$ are homotopic
in $\hat{M}$ to $f$ and $f'$, respectively. Since, $f$ and
$f'$ aren't homotopic in $\hat{M}$ we have a contradiction.
\qed{proof of Lemma \ref{neighbourhood}}

We now prove our main theorem.

\begin{theorem}\label{theoremA}
Let $N$ be a compact, orientable, atoroidal, irreducible 3-manifold
with boundary.  Suppose that $N$ contains an essential, boundary
incompressible annulus whose core curve is not
homotopic into a torus boundary component of $\del N$.  Let $B$ be a
component of the interior of $AH(\pi _1(N))$.  Then there is a
representation $\rho$ in $\overline{B}$ such that for any sufficiently
small neighborhood $V$ of $\rho$ in $AH(\pi _1(N))$ the set $V\cap B$
is disconnected. 
\end{theorem}

\proof

We recall our standing assumption that if $\rho \in B$ then the
marking map $f_\rho:N \longrightarrow M_\rho$ has a homotopy inverse
that is a homeomorphism onto the interior of $N$. If we want to show
self-bumping at a different component $B'$ we find a new manifold
$N'$ homotopy equivalent to $N$ such that $N'$ and $B'$ have the
above property.
With the exceptions of $N$ being irreducible and atoroidal, all of the
topological assumptions we have made depend only on the homotopy type of
$N$. Since a hyperbolic manifold is automatically irreducible and
atoroidal, $N'$ will also be atoroidal, irreducible and contain an essential,
boundary incompressible annulus. In particular if one component
of $MP(\pi_1(N)$ self-bumps then every component
of $MP(\pi_1(N))$ will self-bump.

By Lemma \ref{immersion}, there is a non-trivial shuffle immersion
$f:N \longrightarrow \hat{M} \subset M$ and $f$ lifts to an embedding
$f'$ in the cover $M'$ associated to $f_*(\pi_1(N))$, with $M'$
homeomorphic to $M$.  Let $\hat{M}_\infty$ be a minimally parabolic
structure on $\hat{M}$ which defines a hyperbolic structure
$M'_\infty$ on $M' = M$.  We use $f$ to pull back a hyperbolic
structure $N_\infty$ $N$ and then $f_\infty:N_\infty \longrightarrow
\hat{M}_\infty$ is an isometric immersion and $f'_\infty:N_\infty
\longrightarrow M'_\infty$ is an isometric embedding. The holonomy,
$\rho_\infty(c)$, of $c$ will be parabolic so by an application of the
second Klein-Maskit combination we can find another parabolic $\gamma$
such that the free product of $\rho_\infty(\pi_1(N))$ and $\gamma$ is
a uniformization $\hat{M}'_\infty$ of $\hat{M}$ such that $M'_\infty$
covers $\hat{M}'_\infty$ and $f'_\infty$ descends to an
embedding. Therefore $f$ and $f'$ satisfy the conditions of Lemma
\ref{neighbourhood} which implies the theorem. 
\qed{proof of Theorem \ref{theoremA}}

\begin{cor} \label{noman}
$\overline{B}$ is not a manifold.
\end{cor}

\proof

If $\overline{B}$ is a manifold then Theorem \ref{theoremA}
implies that $\rho_\infty$ is in the interior of
$\overline{B}$, since it cannot be in the boundary.  However, in
\cite{Sullivan}, Sullivan proves that the interior of $\overline{B}$
is $B$. Since $\rho_\infty$ is not in $B$, $\overline{B}$ is not a manifold.
\qed{proof of Corollary \ref{noman}}

In Theorem \ref{theoremA} we characterized when the components
of $MP(\pi_1(N))$ self-bump. To do so we constructed a representation
where this self-bumping occurs. In our next theorem we describe a
sufficient condition for a representation to be a point of
self-bumping. To describe it we will assume some knowledge of Kleinian
groups.

We now allow $N$ to contain more than one copy of $X$. In particular,
assume that there are $m$ disjoint, pairwise embeddings of $(X, \del
X)$ in $(N, \del N)$, labeled, $X_1, \dots, X_m$. As before we assume
that each $\del_1 X_i$ is an essential, boundary incompressible
annulus and that each core curve, $c_i$ is primitive and not homotopic
to a boundary torus. We further assume that the $c_i$ are
homotopically distinct. For each $i$, 
$1 \leq i \leq m$, choose an integer, $n_i \geq 0$. There is then a
shuffle immersion, $s_{n_1,\dots,n_m}$, that wraps $N$ around $c_i$,
$n_i$ times. 

Let $\hat{M} = M-\cC$. If $\hat \rho$ is a minimally
parabolic, geometrically finite uniformization of $\hat{M}$ then the
space of all minimally parabolic hyperbolic structures on $\hat{M}$,
with the same marking, is
$QD(\hat \rho)$, the {\em quasiconformal deformation space} of
$\hat \rho$. The image of $(s_{n_1,\ldots,n_m})_*(\pi_1(N))$ in
$\pi_1(\hat{M})$ defines a Kleinian subgroup $\Gamma$ of
$\hat{\Gamma}=\hat \rho(\pi _1(\hat M))$ that uniformizes $M$, and a
representation $\rho = \hat \rho \circ (s_{n_1,\ldots, n_m})_*$, with
image $\Gamma$. If $\hat \rho '$ is another
representation in $QD(\hat \rho)$ then $\hat \rho ' \circ
(s_{n_1,\ldots, n_m})_*$ is in $QD(\rho)$, the quasiconformal
deformation space of $\rho$. Therefore $(s_{n_1,\ldots,n_m})_*$
defines a map between $QD(\hat\rho)$ and $QD(\rho)$. Our
previous work shows the following:

\begin{theorem}
\label{defbump}
All representations in $QD(\rho)$ in the image of $QD(\hat\rho)$
under $(s_{n_1,\dots,n_m})_*$ are points of self-bumping for $B$ if
$n_i \neq 0$ for some $i$.
\end{theorem}

Note that $\rho$ will not be minimally parabolic, for the $c_i$ will
all be parabolic in $\Gamma=\rho(\pi _1(N))$. Let $c'_i = \{0\} \times
\{1\} \times S^1 \subset \del _0X_i$. The quotient of the domain of
discontinuity for $\Gamma$ will be a conformal structure on $\del N
-\coprod c'_i$.  As the pinched curves in $\del N$
are determined by the embeddings of the $X_i$, if
$s_{n'_1,\dots,n'_m}$ is another shuffle immersion then the image of
$(s_{n'_1,\dots,n'_m})_*$ will be the same quasiconformal deformation
space, $QD(\rho)$. (While these maps have the same range,
$(s_{n_1,\dots,n_m})_*(\hat{\Gamma}) \neq
(s_{n'_1,\dots,n'_m})_*(\hat{\Gamma})$.) On the other hand, each $X_i$
has an involution which swaps the two components of $\del_0 X_i$. By
performing this involution on some (possibly all) of the $X_i$ we get
a new family of shuffle immersions. The bumping representations
associated to these shuffle immersions will then lie in a different
quasi-conformal deformation space.

We also remark that even in the case where $N$ is an $I$-bundle, Theorem
\ref{defbump} is stronger than McMullen's result in
\cite{McMullen:complexeq}. In McMullen's theorem, all the $c'_i$ must
lie in the same component of $\del N$. Here we have no such
restriction.

We close with the following conjecture.

\begin{Conjecture}\label{necessary}
A representation $\rho$ is a point of self-bumping for $B$ if and
only if there is a non-empty collection of curves $\cal C$ (as above)
in $M$, a shuffle immersion $s$ with respect 
to $\cal C$, and a uniformization $\hat \rho$ of $\hat M=M-{\cal C}$ so that
$\rho = \hat \rho\circ s_*$.
\end{Conjecture}

\bibliographystyle{plain}
\bibliography{bib}
\end{document}